\documentclass[12pt,a4paper]{amsart}

\usepackage{soul}
\usepackage{graphicx}
\usepackage[all]{xy}
\usepackage{xspace}

\newcommand{\flqq}{\textquotedblleft}
\newcommand{\frqq}{\textquotedblright\xspace}

\usepackage{DLdef1}

\let\ssec\subsection

\renewcommand {\ssbegin}[1]
  {\refstepcounter{subsection}
  \def \secno {\gdef \secno {}{\ssecfont \thesubsection.\hskip 2ex}%
  }%
   \begin{#1}}

\newcommand{\moplus}{\mathop{\textstyle\bigoplus}\limits}

\newcommand{\fpoi}{{\mathfrak{poi}}}
\newcommand{\wtg}{{\widetilde{\fg}}}

\newcommand{\del}{{\partial}}
\newcommand{\ovR}{{\overline{R}}}

\renewcommand{\Sh}{\mathop{\mathrm{Sh}}\nolimits}

\begin{document}

\title{Shapovalov determinant for loop superalgebras}

\author{Alexei Lebedev${}^1$, Dimitry Leites${}^2$}

\address{${}^1$Equa Simulation AB, Stockholm, Sweden\\
\email{yorool@mail.ru}\\
${}^2$Department of Mathematics, University of Stockholm, Roslagsv.
101, Kr\"aftriket hus 6, SE-106 91 Stockholm, Sweden\\
mleites@math.su.se}

\thanks{We are thankful to MPIMiS, Leipzig, and
the International Max Planck Research School affiliated to it for
financial support and most creative environment.}

\begin{abstract} Let a given finite dimensional simple Lie superalgebra $\fg$
possess an even invariant non-degenerate supersymmetric bilinear
form. We show how to recover the quadratic Casimir element for the
Kac--Moody superalgebra related to the loop superalgebra with values
in $\fg$ from the quadratic Casimir element for $\fg$. Our main tool
here is an explicit Wick normal form of the even quadratic Casimir
operator for the Kac--Moody superalgebra associated with $\fg$; this
Wick normal form is of independent interest.

If $\fg$ possesses an odd invariant supersymmetric bilinear form we
compute the cubic Casimir element.

In addition to the cases of Lie superalgebras $\fg(A)$ with Cartan
matrix $A$ for which the answer was known, we consider the Poisson
Lie superalgebra $\fpoi(0|n)$ and the related Kac--Moody
superalgebra.
\end{abstract}



\maketitle

\markboth{\itshape Alexei Lebedev\textup{,} Dimitry
Leites}{{\itshape Shapovalov determinant for Kac--Moody algebras}}

\thispagestyle{empty}

\section{Introduction}

\subsection{Motivations: physics} For a very lucid detailed
exposition of these motivations in the non-super case, see the
text-book by P.~Di Francesco, P.~Mathieu and D.~S\'en\'echal
\cite{DMS}.

Observe that, in various applications, the simple algebra which
initiated the study is often less interesting than certain its
\lq\lq relatives" (its nontrivial central extension, algebra of
derivations and the result of iterations of these constructions). In
what follows, having the simple object at the center of attention,
we keep its relatives in view.

{\bf 1) Stringy algebras}. In a seminal paper \cite{BPZ},
A.~Belavin, A.~Polyakov and A.~Zamolodchikov observed that the
infinite number of generators of the conformal group in the
two-dimensional case generate the Ward identities for correlation
functions, and these differential equations (Ward identities)
completely specify the behavior of the correlation functions. The
components of the stress-energy operator in the conformal field
theory form, together with the central charge, the Virasoro algebra;
this reduces the study of the conformal theory to the study of
(irreducible) highest weight representations of the Virasoro
algebra.

V. Dotsenko and V. Fateev \cite{DF} explicitly constructed a large
class of conformal theories, so-called {\it minimal models}. In
their study, the complete description of (irreducible) unitarizable
highest-weight representations of a real form of the Virasoro
algebra is vital.

{\bf 2) Loop algebras}. Bosonization of free fermions with spin and
the internal symmetry group $G$ provide with an example of a
nontrivial conformal theory based on the Kac--Moody
algebra\footnote{The formal definition, see \cite{K}, is a bit
complicated: $\widehat{\fg^{(1)}}$ is a certain subalgebra of the
Lie algebra of derivations of the central extension
$\widetilde{\fg^{(1)}}$ of the loop algebra $\fg^{(1)}$ of maps
$S^1\tto \fg$, where $\fg=\Lie(G)$ is a simple finite dimensional
complex Lie algebra, cf. \eqref{bracket}.} $\widehat{\fg^{(1)}}$.
The components of the stress-energy operator for these theories,
built out of quadratic forms of fermion current operators for
$G^{(1)}$, satisfy the relations for the Virasoro algebra with
central charge $C =\frac{\dim \fg}{k + c_\fg}$, where $k$ is the
value of the central charge of $\widehat{\fg^{(1)}}$ and $c_\fg$ is
the value of the (quadratic) Casimir operator of $\fg$ in the
adjoint representation, cf. \cite{GO}. The Hamiltonian of the WZW
model can also be built with quadratic forms of current operators
and thus also represents a nontrivial conformal field theory. The
differential equations for the multi-point correlation functions of
the WZW primary fields are the Knizhnik-Zamolodchikov equations
introduced in \cite{KZ}.

In all these cases the description of the irreducibles is performed
by means of the Shapovalov determinant; in what follows we recall
its definition and ways to compute it.

{\bf 3) Super versions}. For (relatives of) simple stringy
superalgebras, the Shapovalov determinant was computed, except for a
number of cases. (A recent paper \cite{KW} claims to solve all the
cases, but has various omissions.) In this paper, we consider the
Kac-Moody superalgebras; for the stringy Lie superalgebras, see
\cite{GLL}.

\ssec{Means: mathematics} {\bf The Casimir elements}. From the above
we already see that Casimir elements (the elements of the center of
the universal enveloping algebra or its completion) are very
important. In many problems, it suffices to know only the quadratic
Casimir, but we have to have it explicitly.

{\bf The Shapovalov determinant}  is a useful tool for verifying if
certain representations of Lie algebras and Lie superalgebras $\fg$
with vacuum vector (either highest or lowest weight one), namely,
the Verma modules, are irreducible or not, and even for {\it
constructing} certain irreducible modules with vacuum vector.

To \so{define} the Shapovalov determinant, the algebra $\fg$ must be
rather \lq\lq symmetric", i.e.,
\begin{equation}\label{eq1}
\begin{minipage}[c]{.78\textwidth}
$\fg$ must possess a Cartan subalgebra (i.e., a maximal nilpotent
subalgebra coinsiding with its normalizer) $\ft$ whose even part is
commutative and diagonalizes $\fg$, and such that the weight $0$
eigenspace of $\fg$
(relative $\ft_\ev$) coinsides with $\ft$;\\
$\fg$ must possess an involutive (or super-involutive if we respect
the Sign Rule, see below) anti-automorphism $\sigma$ which
interchanges the root vectors (with respect to $\ft_\ev$)
of~opposite sign.
\end{minipage}
\end{equation}
Recall that an {\it anti-automorphism}\index{anti-automorphism} of a
given Lie superalgebra $\fg$ is a linear map $\sigma\colon\fg\tto
\fg$ that for every $x, y\in \fg$, satisfies
\begin{equation}\label{aa}\sigma([x,
y])=\begin{cases}(-1)^{p(x)p(y)}[\sigma(y), \sigma(x)]&\text{if we
respect the Sign Rule},\\ [\sigma(y), \sigma(x)]&\text{if we ignore
the Sign Rule}.\end{cases} \end{equation} (An endomorphism $\sigma$
is said to be {\it involutive} if $\sigma^2=\id$ and {\it
super-involutive} if $\sigma^2=(-1)^{p(x)}\id$. Since the Shapovalov
determinant is defined up to a scalar factor, and the Sign Rule only
affects its sign, we may choose a more convenient definition. In
computations it is usually more convenient to ignore the Sign Rule.)

Additionally, the following finiteness condition should hold:
\begin{equation}\label{eq2}
\begin{minipage}[c]{.78\textwidth}
(a) the root spaces of $\fg$ should be finite dimensional;\\
(b) the number of partitions of any positive weight in $U(\fg^+)$,
see \eqref{fgpm}, into positive roots of $\fg$ should be finite.
\end{minipage}
\end{equation}
Nothing else is needed to \so{define} the Shapovalov determinant,
but to \so{compute} it is much easier in presence of the even
quadratic Casimir element $C_2$ (of the center of $U(\fg)$ or its
completion\footnote{The hat over $U$ means that the elements of
$\widehat U(\fg)$ can be infinite sums of elements of $U(\fg)$.}
$\widehat U(\fg)$) or (which is not quite the same if $\dim
\fg=\infty$, but suffices for our purposes) in presence of the
non-degenerate \so{EVEN} bilinear form $B$ on~$\fg$. \so{In presence
of such an element $C_2$, the Shapovalov determinant is a product of
linear terms}; for various cases where this statement is proved, see
\cite{KK, GL1, Go2, Sh}.

\ssec{Cases where an even quadratic Casimir exists} Kac and
Kazhdan \cite{KK} computed the Shapovalov determinant for any Lie
algebra with symmetrizable Cartan matrix (with arbitrary (complex)
entries); their technique is literally applicable to Lie
superalgebras with symmetrizable Cartan matrix, as mentioned in
\cite{GL1} and expressed in detail in \cite{Go2}.

Moreover, Kac and Kazhdan used the Shapovalov determinant to
describe the {\it Jantzen filtration} of the Verma modules over the
Lie algebras they considered.

The technique of Kac and Kazhdan can be applied even in the absence
of symmetrizable Cartan matrix; its presence only makes computation
of the needed values of the even quadratic Casimir easier. So it is
reasonable to look around for Lie (super)algebras without
symmetrizable Cartan matrix but with a non-degenerate even bilinear
form. Grozman and Leites considered all such simple Lie
superalgebras among $\Zee$-graded of polynomial growth and their
relatives (for the same reasons that Kac--Moody algebras are \lq\lq
better" than simple loop algebras, the finite dimensional Poisson
Lie superalgebras $\fpoi(0|n)$ are \lq\lq better" than their simple
relatives $\fh'(0|n)$, where $\fh(0|n)=\fpoi(0|n)/\mathfrak{center}$
and $\fg'=[\fg, \fg]$): One stringy superalgebra, $\fk^L(1|6)$
(physicists call it $N=6$ Neveu-Schwarz superalgebra), finite
dimensional Poisson superalgebras $\fpoi(0|2n)$ and Kac--Moody
algebras associated with $\fpoi(0|2n)$. More precisely, \cite{GL1}
contains computations of the quadratic Casimirs $C_2$ and
description of the irreducible Verma modules over $\fk^L(1|6)$ and
$\fpoi(0|2n)$.

Gorelik and Serganova \cite{GS1} wrote a sequel to \cite{GL1} having
produced a more explicit expression of the Shapovalov determinant
than that of \cite{GL1}; their rather nontrivial result is an
explicit description of {\it Jantzen's filtration} for the Verma
modules over $\fpoi(0|2n)$.

There are many (perhaps, undescribably many) examples of {\it
filtered} Lie (super)algebras of polynomial growth (i.e., such that
the associated graded Lie (super)algebras are of polynomial growth)
and possessing a $C_2$, cf. \cite{GL0, Ko}. The Shapovalov
determinant is computed only for one (the simplest) of such
algebras, namely, for $\fgl(\lambda)$, where $\lambda \in \Cee$,
 and only for the simplest types of Verma
modules (\cite{Sh}).

\ssec{Cases where no even quadratic Casimirs exist}

{\bf 1) There is an odd non-degenerate bilinear form}. Such are the
queer Lie superalgebra $\fq(n)$ ---  a non-trivial super analog of
$\fgl(n)$, the Poisson superalgebras $\fpoi(0|2n+1)$, and the
Kac-Moody superalgebras associated with them.

In the 1990s, Grozman and Leites conjectured that since $\fq(n)$ and
$\fpoi(0|2n+1)$ are super analogs of $\fgl(n)$, it is possible to
compute their Shapovalov determinant (which is not easy even to
define in these and similar cases), but made a mistake computing it
and decided that the terms it factorizes into can be of any degree.

Gorelik skilfully used the {\it anti-center} and suggested an
elegant proof of the fact that, for $\fq(n)$, the Shapovalov
determinant factorizes in the product of \so{linear} polynomials,
see \cite{Go1}.

Among the simple $\Zee$-graded Lie superalgebras of polynomial
growth that possess a non-degenerate invariant bilinear odd form
there is only one stringy Lie superalgebra (we will consider it
in\cite{GLL}) as well as $\fq(n)$, $\fpoi(0|2n+1)$, and Kac-Moody
algebras associated with them.

{\bf 2) There are no non-degenerate bilinear forms}. Such are most
of the stringy Lie superalgebras (see \cite{GLL}), and Lie
superalgebras $\fq^{(2)}(n)$ with non-symmetrizable Cartan matrices
recently considered in \cite{GS2}.

\ssec{Our results and open problems} Thus, there remained to
consider, first of all, the Lie (super)algebras with properties
\eqref{eq1}, \eqref{eq2}. The Kac--Moody (super)algebras associated
with the loop (super)algebras with values in finite dimensional
\lq\lq symmetric" Lie (super)algebras are among such algebras, so
here we explicitly describe the Wick normal form of the Casimir
operator for the Kac--Moody superalgebras $\widehat{\fg^{(1)}}$ (as
well as \lq\lq twisted" Kac--Moody superalgebras
$\widehat{\fg^{(r)}}$ ) in terms of the the Casimir operator for a
finite dimensional simple Lie superalgebra $\fg$. This implies a
description of the Shapovalov determinant for $\widehat{\fg^{(1)}}$
previously known only for Kac--Moody (super)algebras with Cartan
matrix. We also conjecture a description of the  Shapovalov
determinant for Kac--Moody (super)algebra related with $\fpoi(0|2n)$
for $n>2$.

We will show that if $\fg$ possesses an {\bf odd} invariant
non-degenerate symmetric bilinear form, then $U(\fg)$ (or
$\widehat{U(\fg)}$) contains a {\it cubic} central element. If $\dim
\fg<\infty$, this implies, as a rule, that the Shapovalov
determinant factorizes in a product of factors of degree $\leq 2$.
We also conjecture a description of the  Shapovalov determinant for
$\fpoi(0|2n+1)$.

\section{Background: Lie superalgebras}
For a detailed background, see \cite{GLS}. In what follows the
ground field is $\Cee$.

\ssec{The Poisson superalgebra} Let $G(m)$ be the Grassmann
superalgebra generated by $\theta_1, \dots, \theta_m$. The Poisson
Lie superalgebra (do not confuse with the Poisson-Lie
(super)algebra) $\fpoi(0|m)$ has the same superspace as $G(m)$ and
the (Poisson) bracket is given by
\begin{equation}
\label{theta}
\renewcommand{\arraystretch}{1.4}
\begin{split}
    \{f, g\}_{P.b.}=
    (-1)^{p(f)}\sum\limits_{j\leq m}\ \frac{\partial f}{\partial
\theta_j}\ \frac{\partial g}{\partial \theta_j}\; \text{ for any }\;
f, g\in {\mathbb C} [\theta].
\end{split}
\end{equation}
It is often more convenient to re-denote the $\theta$s and set (over
$\Cee$, over $\Ree$ such a transformation is impossible so the
brackets \eqref{theta} and \eqref{xi} are not equivalent over
$\Ree$)
\begin{equation}
\label{4eq5}
\renewcommand{\arraystretch}{1.4}
\begin{array}{ll}
\begin{cases}\xi_j=\frac{1}{\sqrt{2}}(\theta_{j}-\sqrt{-1}\theta_{r+j})&\cr
\eta_j=\frac{1}{ \sqrt{2}}(\theta_{j}+\sqrt{-1}\theta_{r+j})& \cr
\end{cases}&{\rm  for}\;
j\leq r= \left[\displaystyle\frac{m}{2}\right],\\
\hphantom{\Bigg\{\mbox{}_j}\theta =\theta_{2r+1}& \text{if $m$ is
odd.}
\end{array}
\end{equation}
and accordingly modify the bracket (if $m=2r$, there is no term
with $\theta$):
\begin{equation}
\label{xi}
\renewcommand{\arraystretch}{1.4}
\begin{split}
    \{f, g\}_{P.b.}= 
    (-1)^{p(f)}\bigg(\sum\limits_{j\leq m}\left(\frac{\partial
f}{\partial\xi_j} \frac{\partial g}{\partial
\eta_j}+\frac{\partial f}{\partial \eta_j}\ \frac{\partial
g}{\partial \xi_j}\right)+\frac{\partial f}{\partial \theta}\
\frac{\partial g}{\partial \theta}\bigg).
\end{split}
\end{equation}

The quotient of $\fpoi(0|n)$ modulo center is the Lie superalgebra
$\fh(0|n)$ of hamiltonian vector fields generated by functions:
\begin{equation}
\label{ham}
\renewcommand{\arraystretch}{1.4}
\begin{split}
    H_f:=
    (-1)^{p(f)}\bigg(\sum\limits_{j\leq m}\left(\frac{\partial
f}{\partial\xi_j} \frac{\partial}{\partial \eta_j}+\frac{\partial
f}{\partial \eta_j}\ \frac{\partial }{\partial
\xi_j}\right)+\frac{\partial f}{\partial \theta}\ \frac{\partial
}{\partial \theta}\bigg).
\end{split}
\end{equation}

\ssec{The integral} The still conventional notation
$d^n\theta:=d\theta_1\dots d\theta_n$ for the volume element in the
Berezin integral is totally wrong as was clear already in 1966 from
the explicit form of the Berezinian of the Jacobi matrix of the
coordinate change. A reasonable notation with compulsory indication
of the coordinates is $\vvol(\theta)$; for more motivations, see
\cite{Del, GLS}.

On $\fpoi(0|n)$, or rather on the space of generating functions, the
integral --- equal to the coefficient of the highest term in Taylor
series expansion in $\theta$ --- determines a non-degenerate
invariant supersymmetric bilinear form, of the same parity as $n$,
given by
$$
(f\mid g):=\int fg~\vvol(\theta).
$$

\ssec{Cartan subalgebras, maximal tori, roots and
coroots}\index{Cartan subalgebra}\index{maximal
torus}\index{root}\index{coroot} In \cite{BNO, PS}, it is shown that
the Cartan subalgebras of simple and certain non-simple (like
$\fpoi$, $\fq$ and their subquotients) finite dimensional Lie
superalgebras are conjugate by inner automorphisms. We always fix a
Cartan subalgebra $\ft$ (for example, for $\fpoi(0|2n)$, we take
$\ft=\Cee[\xi_1\eta_1, \dots , \xi_n\eta_n]$). For any
$\alpha\in\ft_\ev^*$, set
 $$
\fg_{\alpha}:=\big\{x\in \fg\mid
\left(\ad_h-\alpha(h)\right)^{N}(x)=0\;\text{ for a sufficiently
large $N$ and every $h\in\ft_\ev^*$}\big\}.
 $$
Then, as is not difficult to see,
\begin{equation} \label{fg=}
\fg=\moplus_{\alpha\in \ft_\ev^*}\fg_{\alpha}\text{~~ and
$\ft\subset\fg_0$}.
\end{equation}

\sssbegin{Remark} In what follows we only consider the algebras for
which $\ft=\fg_0$. \end{Remark}

Denote by $R$ the set of non-zero functionals $\alpha\in\ft_\ev^*$
for which $\dim \fg_{\alpha}\neq 0$; this $R$ is called the set of
{\it roots} of $\fg$. For the algebras we consider,
\[
\text{there exists an $H\in\ft$ such that $\alpha(H)\in\Ree$ for
all $\alpha\in R$ and $\alpha(H)\neq 0$ for $\alpha\neq 0$.}
\]
This allows us to split the roots into positive and negative ones by
setting
\begin{equation} \label{fgpm}
\renewcommand{\arraystretch}{1.4}
\begin{array}{l}
 R^{+}:=\{\alpha\in R\mid \alpha(H)>0\},\quad R^{-}:=\{\alpha\in
R\mid \alpha(H)<0\};\\
 \fg^{\pm}= \moplus_{\alpha\in R^{\pm}}\fg_{\alpha}.
\end{array}
\end{equation}
If $\ft=\ft_\ev$ and commutative, we can identify $U(\ft)$ with
$S(\ft)$; let $HC$ be the {\it Harish-Chandra
projection},\index{Harish-Chandra projection} i.e., the projection
on the first direct summand in the following decomposition:
\begin{equation} \label{HC}
U(\fg)\simeq U(\ft)\moplus \left(\fg^-U(\fg)+U(\fg)\fg^+
\right)\tto S(\ft)\simeq U(\ft).
\end{equation}

\ssec{The Shapovalov determinant}\index{Shapovalov determinant}
\underline{1) $\ft$ is purely even and commutative.} For a fixed
$\lambda\in\ft^*$, set $M^\lambda=U(\fg)/I$, where $I$ is the left
ideal generated by $\fg^+$ and the elements
$$
h-\lambda(h)\;\text{ for any $h\in\ft$}.
$$
The $\fg$-module $M^\lambda$ is called the {\it  Verma
module}\index{Verma module} with highest weight $\lambda$.
Obviously, as spaces, $M^\lambda\simeq U(\fg^-)m^\lambda$, where
$m^\lambda$ is the vacuum vector; and the $U(\fg^-)$-action on
$m^\lambda$ is faithful, i.e., without kernel. The anti-automorphism
$\sigma$, see \eqref{aa}, can be (uniquely) extended to an
anti-automorphism of $U(\fg)$:
\begin{equation}
\label{aaU} \sigma(x_1\ot\dots\ot x_k)=\begin{cases}
(-1)^{\sum\limits_{i<j} p(x_i)p(x_j)}\sigma(x_k)\ot\dots\ot
\sigma(x_1)&\text{minding the Sign Rule;}\\\sigma(x_k)\ot\dots\ot
\sigma(x_1),&\text{ignoring the Sign Rule.}\end{cases}
\end{equation}
On $U(\fg^-)$, we define the following bilinear form with values in
$U(\ft)=S(\ft)$:
\begin{equation}
\label{shdet}(Xm\mid Ym)=HC(a(X)Y).
\end{equation}
Obviously, the distinct weight subspaces of $U(\fg^-)$ are
orthogonal with respect to this inner product $(\cdot\mid\cdot)$
defined by (\ref{shdet}), and $M^\lambda$ is irreducible if and only
if $(\cdot\mid\cdot)$ is non-degenerate. Each determinant $\Sh_\mu$
of the Gram matrix of the restriction of the form $(\cdot\mid\cdot)$
on the subspace $U(\fg^-)(-\chi)$ of weight $-\chi$ is referred to
as the {\it Shapovalov determinant}. Since $S(\ft)$ is commutative
algebra, these determinants are well-defined (as polynomials in
$\ft$); if we do not fix a basis of $U(\fg^-)(-\chi)$ (as we do in
what follows), they are determined up to a scalar factor.

On $M^\lambda$, define the following bilinear form (with values in
$\Cee$):
\begin{equation}\label{shdet-lambda}
(Xm^\lambda\mid Ym^\lambda)_\lambda=(X|Y)(\lambda)\equiv
HC(\sigma(X)Y)(\lambda)\quad\text{for any~}X,Y\in U(\fg^-).
\end{equation}
This definition is natural in the sense that for any $X,X',Y,Y'\in
U(\fg)$ (not only of $U(\fg^-)$), if $Xm^\lambda=X'm^\lambda$ and
$Ym^\lambda=Y'm^\lambda$, then
$$HC(\sigma(X)Y)(\lambda)=HC(\sigma(X')Y')(\lambda).$$ Obviously

1) the elements of distinct weights of $M^\lambda$ are orthogonal to
each other relative $(\cdot\mid\cdot)_\lambda$;

2) the restriction of the bilinear form $(\cdot\mid\cdot)_\lambda$
onto the subspace $M^\lambda(\lambda-\chi)$ of elements of weight
$\lambda-\chi$ is degenerate if and only if $\Sh_\chi(\lambda)=0$.

\sssbegin{Statement} 1) Any non-trivial submodule of $M^\lambda$
contains a non-zero element $v$ whose weight is distinct from
$\lambda$ and $\fg^+v=0$ (such an element $v$ is said to be a {\em
singular} vector of $M^\lambda$).

2) If $v\in M^\lambda$ is a singular vector, then $U(\fg^-)v$
--- is a non-trivial submodule of $M^\lambda$.

3) $x\in M^\lambda$~is an isotropic element of the form
$(\cdot\mid\cdot)_\lambda$ if and only if $x$ can be represented as
$Av$, where $A\in U(\fg^-)$, and $v$ is a singular vector of
$M^\lambda$.

4) If $v\in M^\lambda$ is a singular vector of weight $\mu$ and $C$
is a central element of $U(\fg)$, then
$HC(C)(\lambda)=HC(C)(\mu)$.\end{Statement}

Thanks to this statement to describe all irreducible Verma modules
is the same as to compute all Shapovalov determinants. Moreover, we
see that Casimir elements help to compute Shapovalov determinants.

\underline{2) $\ft_\od\neq 0$.} (This case is considered in more
detail separately.) Among finite dimensional simple Lie
superalgebras, this only happens with $\fpsq(n)$ and $\fh'(0|2n-1)$.
M.~Gorelik \cite{Go1} considered the case of $\fpsq(n)$; that of
$\fh'(0|2n-1)$ might, in theory, be obtained from her results by
contraction, but in applications an {\it explicit} formula is needed
instead of a small talk. In what follows we offer a conjectural
formula for the Shapovalov determinant for $\fpoi(0|2n+1)$;
conjecturally the answer for $\fh'(0|2n+1)$ for $n>1$ is analogous.

\section[Casimir operators: an even invariant form]{Casimir operators
on Lie superalgebras: the case of an even invariant form} In this
section, let $\fg$ be a finite dimensional Lie superalgebra, and
$(\cdot\mid \cdot)$ an invariant even supersymmetric non-degenerate
bilinear form on $\fg$. Let $\{e_i\}_{i=1}^d$ be a basis of $\fg$.
We set
$$
\begin{array}{c}
a_{ij}=(e_i\mid e_j) \text{ and }  (b_{ij})=(a_{ij})^{-1} \text{ the
inverse matrix}.\end{array} $$

\ssbegin{Statement} The quadratic element
\begin{equation}\label{Om0}
\Omega_0=\sum\limits_{i,j} b_{ij}e_i\ot e_j\in U(\fg).
\end{equation}
is a central (Casimir) element of $U(\fg)$, in
particular,
\begin{equation}\label{Om01}
[x,\Omega_0]=0\ \text{ for all }x\in \fg.
\end{equation}
\end{Statement}

\begin{proof} It suffices to prove (\ref{Om01}) for each
$x\in\{e_i\}_{i=1}^d$. Let $c_{ij}^k$ be structure constants in this
basis. Let $p_i=p(e_i)$ be the parities. We have
\begin{equation}
\label{cas}
\begin{split}
{}[e_k,\Omega_0]=&\sum\limits_{i,j} b_{ij}([e_k,e_i]\ot
e_j+(-1)^{p_ip_k}e_i\ot[e_k,e_j])=\\
&\sum\limits_{i,j} b_{ij}\sum\limits_l (c_{ki}^l e_l\ot
e_j+(-1)^{p_ip_k}c_{kj}^l e_i\ot e_l)=\\
&\sum\limits_{i,j} \sum\limits_l
(b_{lj}c_{kl}^i+(-1)^{p_ip_k}b_{il}c_{kl}^j) e_i\ot e_j.
\end{split}
\end{equation}
So, to prove the statement, it suffices to show that
\begin{equation}
\label{b_ij}
 \sum\limits_l
(b_{lj}c_{kl}^i+(-1)^{p_ip_k}b_{il}c_{kl}^j)=0\ \text{ for all }
i,j,k\in\overline{1,d}.
\end{equation}
The invariance of the form $(\cdot\mid \cdot)$ implies that
$$
0=([e_p,e_k]\mid e_r)-(e_p\mid [e_k,e_r])=\sum\limits_s (c_{pk}^s
a_{sr}-c_{kr}^s a_{ps})\ \text{ for all } k,p,r\in\overline{1,d}.
$$
Let us multiply this equality by $b_{ip}b_{rj}$ and sum over $p$
and $r$. Since $(a_{ij})(b_{ij})=1_n$, we get
 $$
0=\sum\limits_p c_{pk}^j b_{ip}-\sum\limits_r c_{kr}^i b_{rj}=
\sum\limits_l (c_{lk}^j b_{il}-c_{kl}^i b_{lj})=-\sum\limits_l
((-1)^{p_kp_l}c_{kl}^j b_{il}+c_{kl}^i b_{lj}).
 $$
Since $(\cdot\mid \cdot)$ is even, the supermatrix $(b_{ij})$ is
also even, and the first term in the last sum can be non-zero only
if $p_i=p_l$, so this equality is equivalent to
(\ref{b_ij}).\end{proof}

Given a loop superalgebra $\fg^{(1)}=\fg\ot \Cee[t,t^{-1}]$, define
the {\it Kac--Moody Lie superalgebra}\index{Kac--Moody Lie
superalgebra} \begin{equation}
\label{KMev}\widehat{\fg^{(1)}}=Span(\fg\ot \Cee[t,t^{-1}], \;u,\;
z), \text{ where $u=t\frac{d}{dt}$ and $z$ are even,} \end{equation}
with the following bracket:
\begin{equation}
\label{bracket}
\renewcommand{\arraystretch}{1.4}
\begin{array}{ll}
 [t^m x,t^n y]=t^{m+n}[x,y]+m\delta_{m,-n}(x\mid y)z &
 \text{for all }x,y\in \fg,\ \ m,n\in\Zee;\\
 {}[z,X]=0 & \text{for all } X\in \widehat{\fg^{(1)}};\\
 {}[u,t^n x]=nt^n x&\text{for all } x\in \fg,\ \ n\in\Zee.
\end{array}
\end{equation}

Set
\begin{equation}
\label{Omega'} \Omega'=\sum\limits_{n=-\infty}^\infty
\sum\limits_{i,j} b_{ij}t^ne_i\ot t^{-n}e_j\in
\widehat{U}(\widehat{\fg^{(1)}}).
\end{equation}

\ssbegin{Statement}\label{Cas0} $[X,2u\ot z+\Omega']=0$ for any
$X\in \widehat{\fg^{(1)}}$.\end{Statement}

\begin{proof} It is easy to check this for $X=u,z$. Now let $X=t^m
e_k$. Let us compute $[X,\Omega']$ up to terms with $z$: we get
\[
\begin{split}
 &\sum\limits_{n=-\infty}^\infty \sum\limits_{i,j}
 b_{ij}(t^{m+n}[e_k,e_i]\ot t^{-n}e_j+
 (-1)^{p_ip_k}t^ne_i\ot t^{m-n}[e_k,e_j])=\\
 &\sum\limits_{n=-\infty}^\infty \sum\limits_{i,j}
 b_{ij}(t^{m+n}[e_k,e_i]\ot t^{-n}e_j+
 (-1)^{p_ip_k}t^{m+n}e_i\ot t^{-n}[e_k,e_j])=\\
 &\sum\limits_{n=-\infty}^\infty \text{\lq\lq}(t^{m+n}\ot t^{-n})\cdot
 [e_k, \Omega_0]\text{\rq\rq}=0,
\end{split}
\]
where the term in quotation marks is understood as follows: having
represented $[e_k, \Omega_0]$ as the sum of quadratic elements, we
multiply, in each of them, the first factor by $t^{m+n}$ and the
second one by $t^{-n}$.

Now, let us compute the terms with $z$ in $[X,\Omega']$. We get
these terms only from the terms in the sum with $n=\pm m$, so we get
\begin{equation}
\label{withz}
\begin{split}
&\sum\limits_{i,j} b_{ij}( m(e_k\mid e_i)z\ot
t^me_j+(-1)^{p_kp_i}mt^m e_i\ot (e_k\mid e_j)z)=\\
&m \sum\limits_{i,j}b_{ij}(a_{ki}z\ot
t^me_j+(-1)^{p_k(p_i+p_j)}t^m e_i\ot a_{kj}z).
\end{split}
\end{equation}
Since $b_{ij}\neq 0$ only for $p_i=p_j$, and $(a_{ij})(b_ij)=1_d$,
it follows that (\ref{withz}) is equal to
\[
m\Big(\sum\limits_{j} \delta_{kj}z\ot t^me_j+\sum\limits_{i}
\delta_{ik} t^me_i\ot z\Big)=m(z\ot X+X\ot z)=2mX\ot Z.
\]
Since $[X,\,2u\ot z]=-2mX\ot z$, we see that $[X,\,2u\ot
z+\Omega']=0$.\end{proof}

Now we want to express the quadratic central element of
$\widehat{U}(\widehat{\fg^{(1)}})$ in the {\it Wick normal
form},\index{Wick normal form} i.e., so that, in every tensor
product, the first factor has a non-positive degree with respect to
$u$, and the second one has a non-negative degree. The Wick normal
form of an operator $\Omega$ will be denoted by $:\Omega:$. Set
\[
\Omega^{pm}=\sum\limits_{n=1}^\infty \sum\limits_{i,j}
b_{ij}t^{-n}e_i\ot t^{n}e_j.
\]
Let us compute $[X,\, \Omega_0+2\Omega^{pm}+2u\ot z]$. The bracket
obviously vanishes for $X=u$ and $z$; so let $X=t^mx$. A computation
similar to the previous one shows that $[X, 2u\ot z]$ cancels with
monomials from $[X,\, \Omega_0+2\Omega^{pm}]$ containing $z$, so we
only need to compute $[X, \Omega_0+2\Omega^{pm}]$ up to elements
with $z$. If $X\in \fg$, then $X$ commutes with $\Omega_0$ (as it
was shown before) and, similarly, with every term in the sum over
$n$ in $\Omega^{pm}$. Consider the case $X=t^mx$, where $m>0$. From
$[X,\,\Omega^{pm}]$ we get
\[
\begin{split}
 &\sum\limits_{n=1}^\infty \sum\limits_{i,j}
 b_{ij}(t^{m-n}[e_k,e_i]\ot t^{n}e_j+(-1)^{p_ip_k}t^{-n}e_i\ot t^{m+n}[e_k,e_j])=\\
 &\sum\limits_{n=1}^\infty \sum\limits_{i,j}
 b_{ij}t^{m-n}[e_k,e_i]\ot t^{n}e_j+\sum\limits_{n=m+1}^\infty
 \sum\limits_{i,j} (-1)^{p_ip_k}t^{m-n}e_i\ot t^{n}[e_k,e_j]=\\
 &\sum\limits_{n=1}^m \sum\limits_{i,j} b_{ij}t^{m-n}[e_k,e_i]\ot
 t^{n}e_j+\sum\limits_{n=m+1}^\infty (t^{m+n}\ot t^{-n})\cdot [e_k,
 \Omega_0]=\\
 &\sum\limits_{n=1}^m \sum\limits_{i,j} b_{ij}t^{m-n}[e_k,e_i]\ot t^{n}e_j.
\end{split}
\]
From $[X,\Omega_0+\Omega^{pm}]$ we get
\[
\begin{split}
 &\sum\limits_{n=0}^\infty \sum\limits_{i,j}
 b_{ij}(t^{m-n}[e_k,e_i]\ot t^{n}e_j+(-1)^{p_ip_k}t^{-n}e_i\ot t^{m+n}[e_k,e_j])=\\
 &\sum\limits_{n=0}^\infty \sum\limits_{i,j} b_{ij}t^{m-n}[e_k,e_i]\ot
 t^{n}e_j+\sum\limits_{n=m}^\infty \sum\limits_{i,j}
 (-1)^{p_ip_k}t^{m-n}e_i\ot t^{n}[e_k,e_j]=\\
 &-\sum\limits_{n=0}^{m-1} \sum\limits_{i,j} b_{ij}(-1)^{p_ip_k}t^{m-n}e_i\ot
 t^{n}[e_k,e_j]+\sum\limits_{n=m}^\infty (t^{m+n}\ot t^{-n})\cdot [e_k, \Omega_0]=\\
 &-\sum\limits_{n=0}^{m-1} \sum\limits_{i,j}
 b_{ij}(-1)^{p_ip_k}t^{m-n}e_i\ot t^{n}[e_k,e_j].
\end{split}
\]
Changing $n\mapsto m-n$, $i\leftrightarrow j$, we get
$$
-\sum\limits_{n=1}^{m} \sum\limits_{i,j}
b_{ji}(-1)^{p_jp_k}t^{n}e_j\ot t^{m-n}[e_k,e_i]=
-\sum\limits_{n=1}^{m} \sum\limits_{i,j}
b_{ij}(-1)^{p_j(p_i+p_k)}t^{n}e_j\ot t^{m-n}[e_k,e_i].
$$
So, for $m>0$, we have
\[
\begin{split}
 &[t^m e_k,\;\Omega_0+2\Omega^{pm}+2u\ot z]=\\
 &\sum\limits_{n=1}^{m}
 \sum\limits_{i,j} b_{ij}(t^{m-n}[e_k,e_i]\ot
 t^{n}e_j-(-1)^{p_j(p_i+p_k)}t^{n}e_j\ot t^{m-n}[e_k,e_i])=\\
 &\sum\limits_{n=1}^{m} \sum\limits_{i,j} t^m
 b_{ij}[[e_k,e_i],e_j]=
 mt^m \sum\limits_{i,j} b_{ij}[[e_k,e_i],e_j].
\end{split}
\]
We similarly get the same result for $m<0$. So we get the following:

\ssbegin{Statement}  If, on $\fg$, the map
$$
A\colon x\mapsto \sum\limits_{i,j} b_{ij}[[x,e_i],e_j]
 $$
is equal to $\lambda\, \id$, where $\id$ is the identity operator,
then

{\em 1)} the element
$$
\Omega:=\Omega_0+2\Omega^{pm}+2u\ot z+\lambda u
$$
is central in $\widehat{U}(\widehat{\fg^{(1)}})$.

{\em 2)} Both $\Omega_0$ and $\Omega$ can be represented in the Wick
normal form.

{\em 3)} The linear terms of $:\Omega_0:$ and $:\Omega:$ differ by
$\lambda u$.\end{Statement}

\sssbegin{Remark} 1) Heading 2) holds since $\Omega_0$ is a finite
sum, and $2\Omega^{pm}+2u\ot z+\lambda u$ is already in the normal
form.

2) Conjecture: If $A$ is not scalar, then no non-zero central
quadratic element of $\widehat{U}(\widehat{\fg^{(1)}})$, if any such
exists, can be expressed in the Wick normal form.\end{Remark}

\ssbegin{Statement} The map $A$ commutes with the
$\fg$-action.\end{Statement}

\begin{proof} According to (\ref{b_ij}), we have
\[
\begin{split}
 [e_k, Ax]=&\sum\limits_{i,j} b_{ij}([[[e_k,x],e_i],e_j]+(-1)^{p_k
 p(x)}[[x,[e_k,e_i]],e_j]+\\
 &(-1)^{p_k(p(x)+p_i)}[[x,e_i],[e_k,e_j]])=\\
 &A[e_k,x]+(-1)^{p_k p(x)}\sum\limits_{i,j}
 b_{ij}([[x,[e_k,e_i]],e_j]+(-1)^{p_kp_i}[[x,e_i],[e_k,e_j]])
 \end{split}
\]
and
\[
\begin{split}
 &\sum\limits_{i,j}
 b_{ij}([[x,[e_k,e_i]],e_j]+(-1)^{p_kp_i}[[x,e_i],[e_k,e_j]])=\\
 &\sum\limits_{i,j,l}b_{ij} (c_{ki}^l[[x,e_l],e_j]+
 (-1)^{p_kp_i}c_{kj}^l[[x,e_i],e_l])=\\
 &\sum\limits_{i,j}\left(\sum\limits_{l}
 (b_{lj}c_{kl}^i+(-1)^{p_kp_i}b_{il}c_{kl}^j)\right)[[x,e_i],e_j]=0.
\end{split}
\]
So we see that
$$
[e_k, Ax]=A[e_k,x]\quad\text{for all~}k\in\overline{1,d}\text{~and~}
X\in \fg. \hfill\qed
$$
\noqed\end{proof}

\sssbegin{Remark} From Schur's lemma we deduce that if $\fg$ is
simple, then $A$ is a scalar operator. For example, if
$\fg=\fsl(m|n)$ and $(X\mid Y)+\tr(XY)$, then $A=2(m-n)$. Note that
if $\fg$ is a direct sum of simple algebras, then $A$ is a direct
sum of the corresponding scalar operators, and hence it may be not a
scalar.
\end{Remark}

\ssbegin{Statement} Let $\fg=\fpoi(0|2n)$. Then $A=0$ if $n>1$ and
not a scalar operator if $n=1$.
\end{Statement}

\begin{proof} If $n=1$, then direct computation shows that
$A$ does not act by 0 on homogeneous elements of degree 1 and 2
whereas $A(1)=0$ since $1$ lies in the center.

Note that if $x,y$ are homogenous polynomials such that $(x\mid
y)\neq 0$, then $\deg x+\deg y=2n$. Thus, for any homogenous
polynomial $X$, if $AX\neq 0$, then $\deg AX=\deg X+2n-4$.

If $n=2$, then $A(\theta_1\theta_2\theta_3\theta_4)=0$: This
element must be either $0$ or a homogenous polynomial of degree
$4$, and the only (up to a scalar factor) such polynomial ---
$\theta_1\theta_2\theta_3\theta_4$
--- does not lie in $\fg'$. Since any
element of the basis can be represented as
$$
\frac{\partial}{\partial\theta_{i_1}}\dots
\frac{\partial}{\partial\theta_{i_k}}\theta_1\theta_2\theta_3\theta_4=
\{\theta_{i_1},\{\dots,\{\theta_{i_k},
\theta_1\theta_2\theta_3\theta_4\}\}\dots\}\;\text{for $k=0$, 1, 2,
3, 4},
$$
where $\{\cdot, \cdot\}_{P.b.}$ is the Poisson bracket and $A$
commutes with the $\fpoi(0|2n)$-action, we see that $A=0$.

If $n\geq 3$, then $AX=0$ for any homogenous polynomial $X$ of
degree $\geq 5$. Since by bracketing with a polynomial of degree $5$
one can get any polynomial of degree $\leq 4$, and since $A$
commutes with the algebra action, we deduce that $A=0$.
\end{proof}

\ssec{Twisted Kac--Moody (super)algebras} Now, let us consider the
case where $\fg$ can be represented,  as a linear space, as
$\widetilde{\fg}_0\oplus \widetilde{\fg}_1$ so that
\[
\renewcommand{\arraystretch}{1.4}
\begin{array}{l}
[\widetilde{\fg}_0,\widetilde{\fg}_0]+
[\widetilde{\fg}_1,\widetilde{\fg}_1]\subset
\widetilde{\fg}_0;\qquad
[\widetilde{\fg}_0,\widetilde{\fg}_1]\subset \widetilde{\fg}_1;\\
(\widetilde{\fg}_0\mid \widetilde{\fg}_1)=0.
\end{array}
\]
In other words, the tilde indicates a $\Zee/2$-grading on $\fg$
({\it a priori} having nothing in common with the parity) respected
by the (non-degenerate, even) inner product $(\cdot | \cdot)$.

In this case, we can define a {\it twisted Kac--Moody Lie
superalgebra}\index{Kac--Moody Lie superalgebra! twisted}
\begin{equation}
\label{KMevtw} \widehat{\fg^{(2)}}=Span(\widetilde{\fg}_0\ot
\Cee[t^2,t^{-2}],\;\; \widetilde{\fg}_1\ot t\Cee[t^2,t^{-2}], \;u,
\;z), \end{equation} where $u=t\frac{d}{dt}$ and $z$ are even, with
the bracket (\ref{bracket}). For the list of simple twisted
Kac-Moody {\it super}algebras, see \cite{FLS}.

Let $A$ still be a scalar operator $\lambda I$. Let
$\{e^0_i\}_{i=1}^{d_0}$ be a basis of $\widetilde{\fg}_0$, let
$\{e^1_i\}_{i=1}^{d_1}$ be a basis of $\widetilde{\fg}_1$; set
$$
  a^0_{ij}=(e^0_i\mid e^0_j);\qquad a^1_{ij}=(e^1_i\mid
  e^1_j),\quad\text{and~} (b^0_{ij})=(a^0_{ij})^{-1},
  (b^1_{ij})=(a^1_{ij})^{-1}.
 $$
Set
\[
\begin{split}
 &\Omega_0'=\sum\limits_{i,j} b^0_{ij}e^0_i\ot e^0_j\in U(\widetilde{\fg}_0);\\
 &\Omega^{\pm}=\sum\limits_{n=1}^\infty \left(\sum\limits_{i,j}
b^0_{ij}t^{-2n}e^0_i\ot t^{2n}e^0_j+\sum\limits_{i,j}
b^1_{ij}t^{-2n+1}e^1_i\ot t^{2n-1}e^1_j\right)\in\widehat{U}(\widehat{\fg^{(2)}});\\
\end{split}
\]
\begin{equation} \label{Omega} \Omega=\Omega_0'+2\Omega^{pm}+2u\ot z+\lambda
u. \end{equation}

Note that $\widetilde{\fg}_0$ is a subalgebra of $\fg$, and
$\Omega_0'$ can be computed for $\widetilde{\fg}_0$ in the same way
as $\Omega_0$ was computed for $\fg$. Then, as in the previous
section, we can prove the following

\ssbegin{Statement} The element $\Omega$ (see (\ref{Omega})) is a
central element in $\widehat{U}(\widehat{\fg^{(2)}})$; its linear
part is equal to the linear part of $\Omega_0'$ plus $\lambda u$,
i.e.,
$$
\text{{\it the linear part of $\Omega$ (not counting $\lambda u$)
is defined by $\widetilde{\fg}_0$}  only. }
$$
\end{Statement}

Similarly, if the algebra $\fg$ has a $\Zee/r$-grading
$\fg=\mathop{\bigoplus}\limits_{s=0}^{r-1} \widetilde{\fg}_s$, we
can construct the algebra
$$
\widehat{\fg^{(r)}}=\Span(\mathop{\bigoplus}\limits_{s=0}^{r-1}\wtg_s\otimes
t^s\Cee[t^r,t^{-r}],\;\;u,\;\;z)
$$
with the bracket (\ref{bracket}). As earlier, set:
\[
\begin{split}
 &e^s_i \text{~--- the basis elements in~} \wtg_s;\\
 &a^s_{ij}=(e^s_i|e^s_j);~ (b_{ij}^s)=(a_{ij}^s)^{-1}; \\
 &\Omega_0'=\sum\limits_{i,j} b^0_{ij}e^0_i\ot e^0_j\in U(\widetilde{\fg}_0);\\
 &\Omega^{pm}=\sum\limits_{s=0}^{r-1}\sum\limits_{n=1}^\infty \sum\limits_{i,j}
 b^s_{ij}t^{-2n+s}e^s_i\ot
 t^{2n-s}e^s_j\in\widehat{U}(\widehat{\fg^{(r)}}).
\end{split}
\]
If $A=\lambda I$, then $\Omega=\Omega_0'+2\Omega^{pm}+2u\ot
z+\lambda u$ is a central element.

For the simple finite dimensional Lie (super)algebras having such a
grading with $r>2$ (i.e., $\fsl(2m+1|2n+1)$ with the automorphism
\flqq minus supertransposition\frqq \footnote{Its order seems to be
equal to $4$, but for all superdimensions, except for $(2m+1|2n+1)$,
its order is equal to $2$ modulo the group of inner automorphisms
(see \cite{Se}, where all outer automorphisms of all finite
dimensional Lie superalgebras are listed).}, $\fg=\fosp(4|2;
\sqrt[3]{1})$ and $\fo(8)$ with the gradings induced by the order 3
automorphisms), {\it these gradings are not compatible with the
weight ones}. Therefore the Cartan subalgebra of the Lie
(super)algebra $\widehat{\fg^{(r)}}$ should be construct not on the
base of the Cartan subalgebra $\fh$ of $\fg$, but on the base of the
Cartan subalgebra of $\tilde\fg_0$, which may have nothing in common
with $\fh$. (For example, if $\fg=\fsl(2m+1|2n+1)$, then
$\tilde\fg_0\simeq\fo(2m+1)\oplus\fo(2n+1)$, and
$\fh\cap\tilde\fg_0=\{0\}$.) Therefore, {\bf although the Casimir
elements of $\fg$ and $\widehat{\fg^{(r)}}$ look alike, their
Shapovalov determinants look completely different}.

\section{Casimir operators on Lie superalgebras: the case of an odd form}
In this section, $\fg$ is a finite dimensional Lie superalgebra on
which there is an invariant {\it odd} supersymmetric (here: this is
the same as just symmetric) non-degenerate bilinear form $(\cdot\mid
\cdot)$. Let, as above, $\{e_i\}_{i=1}^d$ be a basis of $\fg$, and
 $$
\begin{array}{c}
a_{ij}=(e_i\mid e_j), \text{ and }  (b_{ij})=(a_{ij})^{-1}
.\end{array}
 $$

\ssbegin{Statement} The cubic element
$$
C_3=\sum(-1)^{p_l}c_{ij}^kb_{im}b_{jl} e_k\ot e_l\ot e_m
$$
is central in $U(\fg)$.\end{Statement}

Proof is similar as that of Statement \ref{Cas0}.

\ssbegin{Remark} Although the expression of $C_3$ looks less
symmetric than that of $\Omega$ from Statement \ref{Cas0}, it is
possible to show that the coefficient
$$F(k,l,m):=\sum(-1)^{p_l}c_{ij}^kb_{im}b_{jl}$$ of $e_k\ot e_l\ot
e_m$ obeys the Sign Rule applied to $e_k\ot e_l\ot e_m$, relative
permutation of the indices $k,l, m$, i.e.,
$$F(k,l,m)=(-1)^{p_kp_l}F(l,k,m)=(-1)^{p_lp_m}F(k,m,l).$$ Therefore
if $\fg$ is not commutative, then the degree of $C_3$ in $U(\fg)$ is
equal to $3$ (and not less).\end{Remark}

Given a Lie superalgebra $\fg$ with an odd invariant form, we can
construct the Lie superalgebra $\widehat{\fg^{(1)}}=\Span(\fg\ot
\Cee[t,t^{-1}], \;u,\; z)$ with relations (\ref{bracket}), but with
an {\it odd} $z$.

\ssbegin{Conjecture} If $\fg$ is a simple Lie superalgebra, then
$\hat{U}(\widehat{\fg^{(1)}})$ possesses a degree $3$ central
element that can be represented in the Wick normal
form.\end{Conjecture}

\section{Kac--Moody-type Lie superalgebras based on queer Lie superalgebras}

\ssec{What to do if the Cartan subalgebra has odd elements} Here we
consider Shapovalov determinant for Lie superalgebras $\fg$ with
Cartan subalgebra $\ft$ such that $0<n_\od=\dim \ft_\od<\infty$ and
$\ft_\ev$ is commutative. Among $\Zee$-graded simple Lie
superalgebras of polynomial grows (and their relatives) the algebras
with these properties are only the Ramond algebras $R(N)=\fk^M(1|N)$
(see \cite{GLL}), the relatives of $\fq(n)$ considered by Gorelik
\cite{Go1}, and the algebras whose Shapovalov determinant nobody
considered yet: relatives of $\fpoi(0|2n+1)$ and Kac--Moody algebras
associated with $\fq(n)$, $\fpoi(0|2n+1)$.

The Shapovalov  form $(\cdot\mid \cdot)$ defined in (\ref{shdet})
takes values in $U(\ft)$. If $\ft$ is purely even, then
$U(\ft)=S(\ft)$ is a commutative algebra, and the determinant of
$(\cdot\mid\cdot)$ is well-defined up to a scalar factor. If
$\ft_\od\neq 0$, it is not even clear how to define the Shapovalov
determinant.

The simplest way to treat this problem is to
\begin{equation}\label{cond}
\text{consider only Verma modules with 1-dimensional vacuum on
which $\ft_\od$ acts trivially.}
\end{equation}
This restriction (\ref{cond}), however, imposes bounds on possible
weights of the vacuum: Under such restriction $[\ft_\od,\ft_\od]$
also trivially acts on the vacuum. If $[\ft_\od,\ft_\od]=\ft_\ev$,
then {\it the only} 1-dimensional $\ft$-module is, up to parity
change, the trivial one.

M.~Gorelik  \cite{Go1} writes that J.~Bernstein suggested to define
the Shapovalov determinant as follows. First of all, note that
$U(\ft)$ is a Clifford algebra over $S(\ft_\ev)$. For a
1-dimensional $\ft_\ev$-module $\Cee(\lambda)$ of weight
$\lambda\in\ft_\ev^*$, we define the {\it vacuum module} to be
\begin{equation}\label{vac}
U(\ft)\otimes_{U(\ft_\ev)} \Cee(\lambda),
\end{equation}
Clearly, $\sdim U(\ft)\otimes_{U(\ft_\ev)} \Cee(\lambda)=
2^{n_\od-1}|2^{n_\od-1}$.

By the Poincar\'e-Birkhoff-Witt theorem, $U(\ft)$ is a filtered
algebra, the associated graded algebra being
$S(\ft_\ev)\otimes\bigwedge(\ft_\od)$. Thus, there is a natural map
(the composition of the contraction with the Berezin integral)
$\int\colon U(\ft)\tto S(\ft_\ev)\otimes\bigwedge^{\dim
\ft_\od}(\ft_\od)\cong S(\ft_\ev)$ defined up to a scalar factor,
but this suffices for us. So if $\ft_\od\neq 0$, another Shapovalov
form, we call it {\it Bernstein's Shapovalov form}\index{Shapovalov
form! Bernstein}\index{BSh-form, see Shapovalov form! Bernstein}
(BSh-form, for short) can be defined to be:
\begin{equation}\label{shdetnew}
B(\cdot \mid\cdot)=\int (\cdot\mid \cdot).
\end{equation}
This form takes values in the commutative algebra $S(\ft_\ev)$, so
its determinant is well-defined. 

Let us consider separately the case where $n_\od=1$. In this case,
$U(\ft)$ is a commutative (not super-commutative!) superalgebra.
Thus, if we use $U(\ft)\otimes_{U(\ft_\ev)} \Cee(\lambda)$  for the
vacuum, the usual Shapovalov form has  determinant well defined up
to a scalar factor. However, this determinant is equal to the one
Gorelik defined \cite{Go1} up to some power of a non-zero element of
$\ft_\od$.

\ssec{Open problems} 1) We can construct {\it Kac--Moody Lie
superalgebras}\index{Kac--Moody Lie superalgebra} from Lie
superalgebras $\fq(n)$ and $\fpoi(0|2n-1)$ (and their \lq\lq
relatives") using the corresponding odd analog of the invariant
Killing form (so $z$ is odd).

In $\fg=\fq(n)$ and $\fsq(n)$, there is a central element
$E=1_{n|n}$, so $[t^nE, \widehat{\fg^{(1)}}]=\Span(t^nE, z)$. For
any element $x\in U(\widehat{\fg^{(1)}}_-)$, we have
$$
\int HC(x\otimes t^nE)=0
$$
(since $HC(x\otimes t^nE)$ can not contain the maximal product of
odd elements of Cartan subalgebra), i.e., this element lies in the
kernel of the BSh-form. If $\fg=\fpq(n)$ (or, again, $\fsq(n)$),
then the invariant form is also degenerate. Since $\int HC(x\otimes
t^nX)$ does not contain $z$ for any element $X$ from the kernel of
the invariant form and any $x\in U(\widehat{\fg^{(1)}}_-)$, it
follows that $t^nX$ lies in the kernel of the BSh-form.

Thus, in order the BSh-form be non-degenerate in the generic case,
we must consider loops with values in $\fg=\fpsq(n)$; the same
applies to $\fpoi(0|2n-1)$: for non-degeneracy of the BSh-form, we
must consider loops with values in $\fh'(0|2n-1)$. We were unable to
compute the Shapovalov determinant for such loops.

2) The expression of the Shapovalov determinant depends on the
system of positive roots; therefore it is vital to know how to pass
from one system to another one. For finite dimensional Lie algebras
and Kac--Moody algebras the passage is performed by elements of the
Weyl group. For Lie superalgebras such a passage is performed by
means of reflections first introduced by Skornyakov, and
independently by Penkov and Serganova \cite{PS}.

Describe all systems of positive roots (at least, an algorithm: how
to pass from one system to another one and an explicit form of at
least one of them) for $\fg=\fpoi$ or $\fh'$ (perhaps $\fg$
augmented by the grading operator).

\section{Explicit formulas}
\ssec{Lie superalgebras with a symmetrizable Cartan matrix} Let
$\fg=\fg(A, I)$, where $I$ is a set of indices labeling odd coroots
$h_i\in\ft$ dual to the simple roots $\alpha_i$, be a Lie
superalgebra with symmetrizable Cartan matrix $A=(A_{ij})$. Since
$A$ is symmetrizable, there exists an invertible diagonal matrix
$D=\diag(d_1, \dots, d_n)$ such that $B=DA$ is a symmetric matrix.

On the root lattice $\Delta=\Span_\Zee(\alpha_1,\dots, \alpha_n)$,
define:

1) a symmetric bilinear form $(\cdot,\cdot)$ by setting
$$
(\alpha_i,\alpha_j)=B_{ij}\quad \text{for any~}i,j=1,\dots,n;
$$

2) a linear function $F$\footnote{The function $F$ is, essentially,
an analog of the function $(\rho,\cdot)$, where if $\dim\fg=\infty$
(the sums are taken multiplicities counted)
$$\rho:=\frac12\left(\mathop{\sum}\limits_{\alpha\in
R^+_\ev}\alpha-\mathop{\sum}\limits_{\alpha\in
R^+_\od}\alpha\right),$$ where $R_\ev^+$ and $R_\od^+$ are the
subsets of $R^+$ consisting of even and odd roots, respectively,
meaning that the corresponding root vectors are even (resp. odd).}
by setting $$ F(\alpha_i)=\frac12B_{ii}\quad \text{for
any~}i,j=1,\dots,n;
$$

3) for any $\gamma=\sum k_i\alpha_i\in\Delta$, set
$$
h_\gamma=\sum k_id_ih_i
$$
(note that, generally, $h_{\alpha_i}\neq h_i$).

Let $R\subset\Delta$ be a root system and $R^+$ the system of
positive roots. Since every root space is either even or odd, every
root can be endowed with a parity. Let $$ \ovR^+_\ev=\{\alpha\in
R^+\mid p(\alpha)=\ev,~\frac{1}{2}\alpha\not\in R\};\quad
\ovR^+_\od=\{\alpha\in R^+\mid p(\alpha)=\od,~2\alpha\not\in R\}.
$$

Define the {\it partition function},\index{function!
partition}\index{Kostant function} also called {\it Kostant
function} and denoted, when Kostant actively worked in
representation theory, $K$ in his honor, on the set of weights of
$U(\fg_-)$ by the formula $$K(\mu)=\dim U(\fg_-)(\mu).$$

\ssbegin{Theorem}[\cite{M, Go2}] For the $\fg(A, I)$-module
$M^{\chi}$, we have:
\begin{equation}
\label{S}
\renewcommand{\arraystretch}{1.4}
\begin{array}{l} \Sh_{\chi}=\\
\prod\limits_{\alpha\in \ovR^+_{\ev}}\prod\limits_{m\in \Nee}
\left(h_\alpha+F(\alpha)-\frac{m}{2}(\alpha,\alpha)\right)^{K(\chi-m\alpha)}
\prod\limits_{\alpha\in \ovR^+_{\od}}\prod\limits_{m=2k+1\in \Nee}
\left(h_\alpha+F(\alpha)-\frac{m}{2}(\alpha,\alpha)\right)^{K(\chi-m\alpha)}
\end{array}\end{equation}
\end{Theorem}

\sssbegin{Remark} Do not use this expression for $\chi$ which are
not weights of $U(\fg^+)$ (for example, for $\chi=k\alpha$, where
$\alpha$ is an odd simple root such that $2\alpha$ is not a root and
$k>1$): in such cases the formula may produce a non-scalar (hence,
wrong) result.\end{Remark}

\ssec{$\fpo(0|2n+1)$} Let the indeterminates be $\xi_i, \eta_i,
\theta$, where $1\leq i\leq n$, and let the Poisson bracket be of
the form
$$
{}[f,g]_{P.b.}=(-1)^{p(f)}\sum\limits_{i=1}^n \left(\frac{\del
f}{\del \xi_i}\frac{\del g}{\del \eta_i}+\frac{\del f}{\del
\eta_i}\frac{\del g}{\del \xi_i}+\frac{\del f}{\del
\theta}\frac{\del g}{\del \theta}\right).
$$
For the Cartan subalgebra take
$\Cee[\xi_1\eta_1,\dots,\xi_n\eta_n,\theta]$. Let $\epsilon_i$ be
the weight of $\xi_i$, where $i=1,\dots, n$. For any weight
$\gamma=\sum c_i\epsilon_i$, set
$$
h_\gamma:=\sum c_i
\xi_1\eta_1\dots\widehat{\xi_i\eta_i}\dots\xi_n\eta_n.
$$
Set also $h_{max}:=\xi_1\eta_1\dots\xi_n\eta_n$.

\sssbegin{Conjecture} Let $\fg=\fpo(0|2n+1)$, where $n>2$, and let
there be selected one of the systems of simple roots containing all
the $\epsilon_i$. Then $\Sh_\chi$ factorizes in the product of
linear terms of the form $h_\gamma$, where the $\gamma$ are the
roots of $\fg$, such that $\chi-\gamma$ positive or zero weights
$h_{max}$.\end{Conjecture}

For $\fpoi(0|3)$ and $\fpoi(0|5)$, all systems of positive roots are
conjugate, and hence it suffices to compute $\Sh_\chi$ for one
system. For $\fpoi(0|3)$, we select the system consisting of
$\epsilon_1$, and for $\fpoi(0|5)$, the one consisting of
$\epsilon_1,\epsilon_2,\epsilon_1\pm\epsilon_2$.

Set:
\[
\renewcommand{\arraystretch}{1.4}
\begin{array}{ll}
h_1=\xi_1\eta_1;\quad h_0={\bf 1} \text{~(the central element of
$\fpo(0|3)$)}&\text{for $\fpoi(0|3)$};\\
h_{10}=\xi_1\eta_1;\quad h_{01}=\xi_2\eta_2;\quad
h_{11}=\xi_1\eta_1\xi_2\eta_2&\text{for $\fpoi(0|5)$}. \end{array}
\]

\sssbegin{Conjecture} If $\fg=\fpoi(0|3)$, then $\Sh_{k\epsilon_1}$
factorizes into the product of linear factors of the form
$$
h_0\quad\text{and}\quad h_1-m/2,\quad\text{where~}m=1,\dots,k.
$$
\end{Conjecture}

\sssbegin{Conjecture} If $\fg=\fpoi(0|5)$, then, for
$\chi=k\epsilon_1+l\epsilon_2$ (such that $\chi$ is positive if
$k,k+l\geq 0$ and at least one of the numbers $k,l$ is positive),
$\Sh_\chi$ factorizes into the product of linear factors of the form
$$
\begin{tabular}{ll}
$h_{11}$;&\\
$h_{01},$&\text{if $k, k+l\geq 1$;}\\
$h_{10},$&\text{if $k+l\geq 1$;}\\
$h_{01}-h_{10}+m,$&\text{for $m\in\Zee$, $1\leq m\leq k$;}\\
$h_{01}+h_{10}-m,$&\text{for $m\in\Zee$, $1\leq m\leq k, (k+l)/2$.}
\end{tabular}
$$
\end{Conjecture}

\sssbegin{Remark} In the two latter conjectures \lq\lq scalar
summands" are meant to be the scalar terms of $U(\fpoi(0|2n+1))$,
not the central elements of $\fpoi(0|2n+1)$ itself.\end{Remark}

\ssec{$\widehat{\fpo(0|2n)^{(1)}}$}
Let the indeterminants be $\xi_i, \eta_i$, where $1\leq i\leq n$,
and the bracket in $\fpo(0|2n)$ is of the form:
$$
{}[f,g]_{P.b.}=(-1)^{p(f)}\sum\limits_{i=1}^n \left(\frac{\del
f}{\del \xi_i}\frac{\del g}{\del \eta_i}+\frac{\del f}{\del
\eta_i}\frac{\del g}{\del \xi_i}\right).
$$
For the Cartan subalgebra we take
$\Cee[\xi_1\eta_1,\dots,\xi_n\eta_n]\oplus\Span(u,z)$. Let
$\epsilon_i$ be the weight of $\xi_i$, where $i=1,\dots, n$, and
$\epsilon'$ be the weight of $t$. For any weight $\gamma=\sum
c_i\epsilon_i+c'\epsilon'$, we set
$$
h_\gamma:=\sum c_i
\xi_1\eta_1\dots\widehat{\xi_i\eta_i}\dots\xi_n\eta_n+c'z.
$$
Set also $h_{max}:=\xi_1\eta_1\dots\xi_n\eta_n$.

\sssbegin{Conjecture} Let $\fg=\widehat{\fpo(0|2n)^{(1)}}$, where
$n>2$, and let there be selected one of the system of positive roots
containing all roots $\epsilon_i, \epsilon'$. Then $\Sh_\chi$
factorizes into the product of linear factors of the form $h_{max}$
and $h_\gamma$, where the $\gamma$ are positive roots of $\fg$ such
that $\chi-\gamma$ are either positive or zero
weights.\end{Conjecture}

\end{document}